\documentclass[11pt]{amsart}
\usepackage{amssymb}
\usepackage[latin1]{inputenc}
\usepackage{graphics}

\usepackage[english,francais]{babel}

\newtheorem{theorem}{Theorem}[section]
\newtheorem{lemma}[theorem]{Lemma}
\newtheorem{e-proposition}[theorem]{Proposition}
\newtheorem{corollary}[theorem]{Corollary}
\newtheorem{e-definition}[theorem]{Definition\rm}
\newtheorem{remark}{\it Remark\/}

\newtheorem{theoreme}{Th\'eor\`eme}[section]

\newtheorem{proposition}[theoreme]{Proposition}

\newcommand{\be}{\begin{equation}}
\newcommand{\ee}{\end{equation}}

\newcommand{\bbR}{{\mathbb R}}

\def\og{\leavevmode\raise.3ex\hbox{$\scriptscriptstyle\langle\!\langle$~}}
\def\fg{\leavevmode\raise.3ex\hbox{~$\!\scriptscriptstyle\,\rangle\!\rangle$}}

\begin{document}
% place in the next line the header (rubrique) chosen for your article,
% if you know it (you can also have 2, format : Header1/Header2
\centerline{}
%\begin{frontmatter}


% \author{Name\thanksref{label2}}
\selectlanguage{english}
\title{Trapping Rossby waves}

% The [label1] can be suppressed if there is only one address for all authors

\selectlanguage{english}
\author[Cheverry]{Christophe Cheverry}
\email{christophe.cheverry@univ-rennes1.fr}
\author[Gallagher]{Isabelle Gallagher}
\email{gallagher@math.jussieu.fr}
\author[Paul]{Thierry Paul}
\email{thierry.paul@ens.fr}
\author[Saint-Raymond]{Laure Saint-Raymond}
\email{Laure.Saint-Raymond@ens.fr}

\address[Cheverry]{Institut Math\'ematique de Rennes, 
Campus de Beaulieu, 263 avenue du G\'en\'eral Leclerc CS 74205 
35042 Rennes Cedex}
\address[Gallagher]{Institut de Math\'ematiques de Jussieu et Universit\'e Paris VII, Case 7012, 2 place Jussieu 
75251 Paris Cedex 05}
\address[Paul]{CNRS and DMA \'Ecole Normale Sup\'erieure, 45 rue d'Ulm, 75230
Paris Cedex 05}
\address[Saint-Raymond]{Universit\'e Paris VI and DMA \'Ecole Normale Sup\'erieure, 45 rue d'Ulm, 75230
Paris Cedex 05}
% If you know the dates of reception, and acceptation you can put them now;
%  idem the name of the person presenting the Note

\maketitle

\begin{abstract}
%\selectlanguage{english}
% Text of abstract in English
Waves associated to large scale oceanic motions are gravity waves (Poincar\'e waves which disperse fast)   and quasigeostrophic waves (Rossby waves). 
In this Note, we show by
semiclassical arguments, that Rossby waves can be trapped and we characterize the corresponding initial conditions.
%{\it To cite this article: C.Cheverry et al. , C. R. Acad. Sci. Paris, Ser. I 340 (2005).}

\vskip 0.5\baselineskip

\end{abstract}
%\end{frontmatter}
% now the Version franÁaise abr\'eg\'ee, if it exists

\selectlanguage{english}
% main text
\section{Introduction and results}
\label{}
Large scale oceanic motions can be described - at first sight - by the linearized Saint-Venant equations for thin layers, with   Coriolis force, 
i.e. by the following two-dimensional system~ \cite{P1}:
\be\label{syst}
\partial_tU+\left(\begin{array}{ccc}
0&\partial_{x_1}&\partial_{x_2}\\
\partial_{x_1}&0&-B(x_2)\\
\partial_{x_2}&B(x_2)&0\end{array}\right)U=0,\ \mbox{where}\  U=U(t,x_1,x_2)=\left(\begin{array}{c}\rho\\u_1\\u_2\end{array}\right)
\ee

\noindent where $B$ is the local vertical component of the Earth rotation vector (depending only on the latitude $x_2$),~$\rho$ denotes 
the fluctuation of water height and $u$  the horizontal velocity field. For the sake of simplicity, we assume that $(x_1,x_2) \in \bbR\times\mathbb T$, 
meaning that we neglect the influence of the lateral boundaries.

Waves associated to that linear system are usually classified in two ~\\families~: First, Poincar\'e waves, which are fast dispersive 
gravity waves; Second, Rossby waves, due to the inhomogeneities of $B$,  which propagate much slower  \cite{DU}, \cite{GA}.

In this Note, we are interested in the propagation of Rossby waves. We show, in agreement with physical observations \cite{P2}, that these waves 
can be trapped in some regions, called ventilation zones, which are not influenced by external dynamics and sources such as continental 
recirculation for instance.

%It is well known that in the $\beta$-plan approximation ($b(x)=x$) there is no trapping of Rossby waves.

%In this note we want to show that, if the $\beta$-plan condition is released ($B(x)$ not linear), there is a (phase-space) trapping region in the limit of large values of $B$.

More precisely, in the limit of large values of $B=b/\epsilon$, we   construct  a co-dimension $1$ submanifold~$\Lambda$ of~$T^*(\bbR\times\mathbb T)$. 
This set contains the $\epsilon$-microlocalization
region inside which, essentially, an initial condition of (\ref{syst})  remains trapped. 
Let us first recall that the $\epsilon$-frequency set of a function $u$ \cite{MA} is the (closed) subset of~$T^*(\bbR\times\mathbb T)$, complement of the set of points $(\underline x,\underline\xi)$ such that there exists a $C^\infty$ function $\chi$, with $\chi(\underline x,
\underline \xi)=1$, such that
\[
\left\Vert\int\chi(\frac{x+y}2,\xi)e^{i\frac{\xi(x-y)}\epsilon}u(y)dyd\xi\right\Vert_{L^2}=O(\epsilon^\infty).
\]

\begin{theorem}\label{main}
Let us consider the system 
\be\label{systempsilon}
\partial_tU+\left(\begin{array}{ccc}
0&\partial_{x_1}&\partial_{x_2}\\
\partial_{x_1}&0&-\frac{b(x_2)}\epsilon\\
\partial_{x_2}&\frac{b(x_2)}\epsilon&0\end{array}\right)U=0
\ \mbox{with an initial condition}\ 
\ee
\be
U|_{t=0}=\left(\begin{array}{c}
\rho^0\\u_1^0\\u_2^0\end{array}\right)
\ee
with $L^2(\bbR\times \mathbb T,dx_1dx_2)$ conditions (periodic in $x_2$). Let $\Lambda=\{F(\xi_1,x_2,\xi_2)=0\}\subset T^*\bbR\times \mathbb T$, where $F$ is defined in Lemma
\ref{pm}.

We    suppose that the $\epsilon$-frequency set of  $U|_{t=0}$ is contained in a compact set $\mathcal C$ satisfying:
\be\label{cond}
\mathcal C\cap\{\xi_1=0\} =\mathcal C\cap\{\xi_2^2+b(x_2)^2=0\}=\emptyset.
\ee
Let us fix a compact set $\Omega$ of $\bbR\times\mathbb T$. Then there exists ($\epsilon$-)pseudo-differential operators $P_\rho^0, P_1^0, P_2^0$ of principal symbols
%\be\label{symb}\scriptsize
$p_\rho^0=ib(x_2)\xi_1(\xi_2^2+\xi_1^2+b^2(x_2))^{-1}$, $p_1^0=-\xi_1\xi_2(\xi_2^2+\xi_1^2+b^2(x_2))^{-1}$, $p_2^0=\xi_1^2(\xi_2^2+\xi_1^2+b^2(x_2))^{-1}$, such that:

%$p_\rho^0=-2i\frac{\sqrt{\xi_1^2+\xi_2^2+b^2(x_2)}b^2(x_2)}{\xi_2^2+b^2(x_2)},\ p_1^0=2\frac{\sqrt{\xi_1^2+\xi_2^2+b^2(x_2)}\xi_1}{\xi_2^2+b^2(x_2)},
%\ p_2^0=2\frac{\sqrt{\xi_1^2+\xi_2^2+b^2(x_2)}\xi_1(\xi_1^2-b^2(x_2))}{\xi_2^2+b^2(x_2)}
%$%\ee,

\begin{enumerate}
\item if the $\epsilon$-frequency set  of $P_\rho^0\rho^0+ P_1^0u_1^0+ P_2^0u_2^0$ does not intersect $\Lambda\cap T^*\Omega$, then  $\exists
T>0$ such that:
\[
||U(\frac t{\epsilon})||_{L^2(\Omega)}=O(\epsilon^\infty)\ \mbox{for}\ t\geq  T .
\]

\item  if the $\epsilon$-frequency set of $P_\rho\rho^0+ P_1u_1^0+ P_2u_2^0$ does  intersect $\Lambda\cap T^*\Omega$, then, $\forall  t\geq 0 $,
$
||U(\frac t{\epsilon})||_{L^2(\Omega)}\neq O(\epsilon^\infty)
$
(in other words  the frequency set of $U(\frac t{\epsilon})$ intersects $T^*\Omega$).
\end{enumerate}

In the case of a WKB initial condition the conclusion is more precise. Consider
\[
U|_{t=0}=\left(\begin{array}{c}
R^0(x)\\U_1^0(x)\\U_2^0(x)\end{array}\right)e^{i\frac{S(x)}\epsilon}.
\]
Suppose than  the (Lagrangian) manifold \\$\{(x,\nabla S(x)),(p_\rho^0 R^0+p_1^0U_1^0+p_2^0U_2^0)(x,\nabla S(x)) \neq 0\}$
intersects $\Lambda\cap T^*\Omega$, then:

\be\label{wkb}
||U(\frac t{\epsilon})||_{L^2(\Omega)}
\sim C(t)+O(\epsilon),\ C(t)>0.
\ee

\end{theorem}

\section{Reduction to a scalar situation}\label{reduction}
Performing first a Fourier analysis in $x_1$ (as the system does not contain explicitly this variable) and looking secondly at modes of the system (\ref{systempsilon}), we can
prove the following Proposition, heart of our results:
\begin{proposition}\label{reduc}
There exist three pseudo-differential operators $T_\pm,\ T_0$, of leading symbols \\ $\tau_{\pm}= \pm\sqrt{\xi_1^2+\xi_2^2+b^2(x_2)},\
\tau_0=\epsilon\frac{b'(x_2)\xi_1}{\xi_1^2+\xi_2^2+b^2(x_2)}$, roots of (\ref{the}), such that, if $u_{2,j}^0$ is microlocalized outside $\xi_1^2-\tau_j^2=0$, we have, $\mbox{for}\  j\in\{-,0,+\}$:
$$\epsilon\partial_tu_{2,j}=iT_ju_{2,j}   \Longrightarrow 
U_j:=\left(\begin{array}{c}(i \epsilon\partial_{x_2}T_j+ \epsilon\partial_{x_1}b(x_2))(\epsilon^2\partial_{x_1}^2+T_j^2)^{-1}\\
-(\epsilon^2\partial_{x_1}\partial_{x_2}+ib(x_2)T_j)(\epsilon^2\partial_{x_1}^2+T_j^2)^{-1}\\
\mathbb Id\end{array}\right)u_{2,j}$$ $$\mbox{ satisfies (\ref{systempsilon}) up to } O(\epsilon^\infty).
$$
\end{proposition}
The proof is a consequence of the result contained in the Appendix.

The following result shows that any initial condition of (\ref{systempsilon}) can be decomposed on the three modes of the last section.
\begin{proposition}
$\forall \rho,u_1, u_2,\ \exists  u_{2,j},\ j\in\{+,0,-\}$ such that:
\be\label{decomp}
\left(\begin{array}{c}\rho\\u_1\\u_2\end{array}\right)=\sum_{j} \left(\begin{array}{c}(i \epsilon\partial_{x_2}T_j+ \epsilon\partial_{x_1}b(x_2))(\epsilon^2\partial_{x_1}^2+T_j^2)^{-1}\\
-(\epsilon^2\partial_{x_1} \partial_{x_2}+ib(x_2)T_j)(\epsilon^2\partial_{x_1}^2+T_j^2)^{-1}\\
\mathbb Id\end{array}\right)u_{2,j}+O(\epsilon^\infty)
\ee
\[=:\sum_{j} {\mathbb Q}^j u_{2,j}+O(\epsilon^\infty) .
\]
\end{proposition}
To prove the Proposition one has just to invert the matrix:
$(
{\mathbb Q}^- \, {\mathbb Q}^0 \, {\mathbb Q}^+ ).$
Semiclassically it is enough to show that the matrix of the leading order symbol is invertible
\be\label{matrix}\scriptsize
\left(\begin{array}{ccc}
{\xi_2\sqrt{\xi_1^2+\xi_2^2+b(x_2)^2} + i\xi_1b(x_2)\over \xi_2^2+b(x_2)^2}&-{ i b(x_2)\over \xi_1}
&{-\xi_2\sqrt{\xi_1^2+\xi_2^2+b(x_2)^2} + i\xi_1b(x_2)\over \xi_2^2+b(x_2)^2}\\
{\xi_1\xi_2+ib(x_2) \sqrt{\xi_1^2+\xi_2^2+b(x_2)^2} \over \xi_2^2+b(x_2)^2}&-{\xi_2 \over \xi_1}
&{\xi_1\xi_2-ib(x_2) \sqrt{\xi_1^2+\xi_2^2+b(x_2)^2} \over \xi_2^2+b(x_2)^2}\\1&1&1
\end{array}\right).
\ee
A simple computation shows that the  jacobian
$
J=\frac{2(\xi_1^2+\xi_2^2+b^2(x_2))^{3/2}}{(\xi_2^2+b^2(x_2))|\xi_1|}\geq 2$.
In particular the inversion of the matrix $(
{\mathbb Q}^- \, {\mathbb Q}^0 \, {\mathbb Q}^+ )$ can be done symbolically at any order and gives the leading order symbols of the operators $\mathbb
P^j:=(P_\rho^j,P_1^j,P_2^j)$
 such
that $u_{2,j}=P_\rho^j\rho+P_1^ju_1+P_2^ju_2,\ j\in\{-,0,+\}$. One gets~:
%\be\scriptsize\label{p
$p_\rho^0=ib(x_2)\xi_1(\xi_2^2+\xi_1^2+b^2(x_2))^{-1},\ p_1^0=-\xi_1\xi_2(\xi_2^2+\xi_1^2+b^2(x_2))^{-1},\  p_2^0=\xi_1^2(\xi_2^2+\xi_1^2+b^2(x_2))^{-1}$.
%$p_\rho^0=-2i\frac{\sqrt{\xi_1^2+\xi_2^2+b^2(x_2)}b^2(x_2)}{\xi_2^2+b^2(x_2)},\ p_1^0=2\frac{\sqrt{\xi_1^2+\xi_2^2+b^2(x_2)}\xi_1}{\xi_2^2+b^2(x_2)},
%\ p_2^0=2\frac{\sqrt{\xi_1^2+\xi_2^2+b^2(x_2)}\xi_1(\xi_1^2-b^2(x_2))}{\xi_2^2+b^2(x_2)}.
%$%\ee

\section{Propagation under the Rossby Hamiltonian $\tau_0$}\label{lambda}
Thanks to  Proposition \ref{reduc} it is enough, as far as  the Rossby mode is concerned, to look at propagation with respect to the Hamiltonian $\tau_0/\epsilon$.
This Hamiltonian being  independent of $x_1$, $\xi_1$ will be conserved. The flow is periodic in the variables $x_2,\xi_2$ (one degree of freedom). 
Therefore, since 
$\dot{x}_1=\frac{b'(x_2)(\xi_2^2-\xi_1^2+b^2(x_2))}{(\xi_2^2+\xi_1^2+b^2(x_2))^2}\mbox{ is periodic (the case of infinite and zero period is treated}$
 in \cite{CGPS}),$\ x_1(t)$ 
will contain a part, linear in time except if
$$F(\xi_1,x_2(0),\xi_2(0)):=\int_0^{period}\frac{b'(x_2(t))(\xi_2(t)^2-\xi_1^2+b^2(x_2(t)))}{(\xi_2(t)^2+\xi_1^2+b^2(x_2(t)))^2}dt=0.$$
%\begin{lemma}\label{pm}
%Define 
%$
%E(\xi_1,x_2,\xi_2)=\frac{b'(x_2)\xi_1}{\xi_2^2+\xi_1^2+b^2(x_2)}.
%$
%Then
%$\displaystyle
%\left\vert F(\xi_1,x_2,\xi_2)\right\vert=\left\vert\int_{x_-}^{x_+}
%\frac{\frac{b'(x)}{E(\xi_1,x_2,\xi_2)}-2\xi_1}{\sqrt{\frac{b'(x)\xi_1}{E(\xi_1,x_2,\xi_2)}-\xi_1^2-b^2(x)}}dx\right\vert,
%$
%where $x_\pm$ satisfy $\frac{b'(x_\pm)\xi_1}{E(\xi_1,x_2,\xi_2)}-\xi_1^2-b^2(x_\pm)=0$, $x_-< x_2< x_+$ and $\frac{b'(x)\xi_1}{E(\xi_1,x_2,\xi_2)}-\xi_1^2-b^2(x)\geq0,\
%x\in[x_-,x_+]$.
%\end{lemma}
%We easily check that $\forall(x_2,\xi_2),\ \pm b(x_2)F(\xi_1,x_2,\xi_2)>0\mbox{ as } \xi_1\to\pm\infty$.
\begin{lemma}\label{pm}
As $b'(x_2)\neq 0,\ b'(x_2)F(\xi_1,x_2,\xi_2)>0\mbox{ as } \xi_1\to\pm\infty,\ <0\mbox{ as }$\\$ \xi_1\to 0$ and is invariant under the flow of  $\frac{\tau_0}\epsilon$.

Define 
$
E(\xi_1,x_2,\xi_2)=\frac{b'(x_2)\xi_1}{\xi_2^2+\xi_1^2+b^2(x_2)}\ \ .
$

Then
$\displaystyle
\left\vert F(\xi_1,x_2,\xi_2)\right\vert=\left\vert\int_{x_-}^{x_+}
\frac{\frac{b'(x)}{E(\xi_1,x_2,\xi_2)}-2\xi_1}{\sqrt{\frac{b'(x)\xi_1}{E(\xi_1,x_2,\xi_2)}-\xi_1^2-b^2(x)}}dx\right\vert,
$
where
%~$x_\pm$ are defined so that
~$]x_-,x_+[$
is the largest interval of~${\mathbb T}$ containing $x_2$ in which~$\frac{b'(x)\xi_1}{E(\xi_1,x_2,\xi_2)}-\xi_1^2-b^2(x)>0$.
\end{lemma}
We define  
$\Lambda:=\{(x_1,x_2;\xi_1,\xi_2)/\ F(\xi_1,x_2,\xi_2)=0\}.$

$\mbox{ Thanks to Lemma \ref{pm} } \Lambda\neq\emptyset\mbox{ and }dim\Lambda=3$.
%Moreover, thanks to the invariance of $F$ under the flow of Hamiltonian $\tau_0$, $\lambda$ is of co-dimension $1$.
\begin{corollary}
Suppose~$b'(x_2) \neq 0$. Then $\vert F(\xi_1, x_2,\xi_2)\vert  \geq\frac{C(x_2,\xi_2)}{\xi_1}$ as $\xi_1\to 0$, with~$C(x_2,\xi_2)>0$. 
This implies that the trapping phenomenon will take place only with initial conditions oscillating enough in $x_1$.
\end{corollary}

\medskip
\begin{remark}: Since the Hamiltonian $E$ does not depend on $x_1$, we can  express it, for each value of $\xi_1$, on the action variable 
$A$: $E(\xi_1,x_2,\xi_2)=H(A,\xi_1)$. This allows to define the function
  $A(\xi_1,x_2,\xi_2)$ by $E(\xi_1,x_2,\xi_2)=H(A(\xi_1,x_2,\xi_2),\xi_1)$. One can easily show that:
\[F(\xi_1,x_2,\xi_2)=\frac{\partial_{\xi_1}H(A,\xi_1)}{\partial_{A}H(A,\xi_1)}|_{A=A(\xi_1,x_2,\xi_2)},\]
  and the following variational characterization of $\Lambda$:
Let us fix the energy to $E$ and let $\Gamma(\xi_1,E)$ be the energy shell $\{E(\xi_1,x_2,\xi_2)=E\}$. Then

\[\Lambda=\bigcup_{E,i}\{(\Gamma(\xi_1^i,E),\xi_1^i\},\]\[
\mbox{ where $\xi_1^i$ is such that the area inside $\Gamma(\xi_1^i,E)$ is extremal.}\]
\end{remark}

\section{Dispersion of Poincar\'e waves and proof of the Theorem}
The proof of the Theorem  involves, for the Rossby modes ($j=0$), the standard result of propagation of the frequency set. If
the initial  frequency set  is such that part of the Rossby mode is  trapped, in particular  if it intersects $\{ b'(x_2)=0\}$, this concludes the proof.

If not we have to prove some  dispersion for the Poincar\'e modes ($j=\pm$) for  times of the order $\frac 1 \epsilon$, for which the theorem of propagation of the frequency set is not enough. But, the system 
being integrable, 
we can perform an expansion on (Bohr-Sommerfeld) eigenvalues of the Hamiltonians  $T_\pm$  and a decomposition of~$U|_{t=0}$ on a compact set of
coherent states (thanks to the condition on the microlocalization of the initial datum) \cite{PA}.

Since  $\tau_\pm=\pm\sqrt{\xi_2^2+b^2(x_2)+\xi_1^2 }+O(\epsilon)$ and we are microlocalized far away from $\xi_2^2+b^2(x_2)+\xi_1^2=0 $, 
we can find  pseudo-differential operators $H_{2\pm}$ of principal symbols 
$\xi_2^2+b^2(x_2)$ such that $T_\pm=\pm\sqrt{H_{2\pm}+\xi_1^2 }$. The Bohr-Sommerfeld
quantization condition (with subsymbol) gives that the eigenvalues of $H_{2\pm}$ are of the form:
\[
\lambda_\pm^k=\lambda_\pm\left((k+\frac 1 2 )\epsilon\right)+\epsilon\mu^{k}_\pm(\xi_1)+O(\epsilon^2),
\]
where $\lambda_\pm$ is the energy $\xi_2^2+b^2(x_2)$ defined on action variable, and $\epsilon\mu^{k}_\pm(\xi_1)\in C^\infty$ is the correction due to the subsymbol.
Propagating at time $t= s /\epsilon$ a function, product of a coherent state at $(q,p)$ (in $x_1$) and  
an eigenfunction of $T_\pm$ (in $x_2 $),  gives rise to expressions of the type:
\[
\int \exp{i\frac{\epsilon(x_1-q)\xi_1\pm(\lambda_\pm^k+\xi^2_1)^{\frac 1 2}s+i\epsilon(\xi_1-p)^2}{\epsilon}}d\xi_1.
\]
The stationary phase lemma then  gives that this integral is $O(\epsilon^\infty)$ except if there exists a stationary point,  given by the conditions:
\[
\xi_1=p\ \mbox{and}\  \epsilon(x_1-q)\pm\frac{(2\xi_1 +\epsilon\partial_{\xi_1}\mu^{k}_\pm)}{2\sqrt{\lambda_\pm^k+\xi^2_1}}s=0 . \]
The second condition gives:
$
2\sqrt{\lambda_\pm^k+\xi^2_1}(x_1-q)=\mp\left(\frac{2p}\epsilon+\partial_{\xi_1}\mu^{k}_\pm\right)s.
$
Therefore, since $p\neq 0$ and the $\lambda^k$\ 's are bounded by the above condition (\ref{cond}) on $\mathcal C$, there is no critical point for $x_1$ in a compact set.

% etc, etc

% \appendix

% \section{}
% \label{}

% The Acknowledgements are an un-numbered section
%\section*{Acknowledgements}
% Acknowledgements text here
\appendix
\section{A microlocal analysis Lemma} In this Appendix we give the crucial lemma for the proof of Proposition \ref{reduc}. It tells us that the principal symbols of 
$T_\pm,\ T_0$ can be computed by solving the symbolic equation associated to 
$$
\left(\epsilon^2\partial_{x_2}^2-b^2(x_2)-i\epsilon^2\frac {b'(x_2)\partial_{x_1}}\tau+\tau^2+\epsilon^2\partial_{x_1}^2\right)u_2=0 \ $$
%\mbox{obtained algebraically from 
%(\ref{syst}), that is, }$$
obtained by (exact) algebraic computations from 
(\ref{syst}), that is
\be\label{the}\tau^2-\xi_1^2-\xi_2^2-b^2+\epsilon\frac{b'\xi_1}\tau=0 .
\ee
%obtained algebraically from (\ref{syst}).
\begin{lemma}
Let $h=h(x,\xi,\sigma)$ be a smooth function such that $\partial_\sigma h|_{h=0}\neq 0$, and let $\tau=\tau(x,\xi)$ be any continuous root of $h(x, \xi,\tau(x,\xi))=0$.

Then there exists a pseudo-differential operator $T=T(x,-i\epsilon\partial_x)$ with principal symbol $\tau(x,\xi)$ such that:
\be
\label{theother}
T\psi=\lambda\psi\ \Longrightarrow\ H(\lambda)\psi=O(\epsilon^\infty)
\ee
where $H(\lambda)$ is a pseudo-differential operator of full symbol $h(x,\xi,\lambda)$. 
\end{lemma}
\textit{Sketch of the proof:} the proof uses essentially   pseudo-differential functional calculus \cite{MA}. 

For any $(x,\xi)$ the principal symbol of $h(x,\xi,T)$ at $(y,\xi')$ is $h(x,\xi,\tau(y,\xi'))$. We can write at first order:
$
H\psi:=\int e^{i\frac{\xi(x-y)}\epsilon}h(x,\xi,T(y,-i\epsilon\partial_y))\psi(y)\frac{d\xi dy}\epsilon=
\int e^{i\frac{\xi(x-y)}\epsilon}e^{i\frac{\xi'(y-y')}\epsilon}h(x,\xi,\tau(y',\xi'))\psi(y')\frac{d\xi d\xi'dy}\epsilon.
$
The integral over $y$ gives $\delta(\xi-\xi')$. Therefore:
$
H\psi(x)=\int e^{i\frac{\xi(x-y)}\epsilon}h(x,\xi,\tau(y,\xi))\psi(y)\frac{d\xi dy}\epsilon.$
So the principal symbol of $H$ is $h(x,\xi,\tau(x,\xi))$ which, by assumption, is $0$.

For the $\epsilon^\infty$ result, it is enough to repeat the same argument taking into account lower order terms $\epsilon^k\tau_k$ for and adding to any symbol a term of the form
$\sum_{k\geq 1}\epsilon^kP_k\tau$ where the $P_k$s are differential operators. We end up with an equation for $\tau_\epsilon\sim \tau+\sum \epsilon^k\tau_k$ of the form:

$
h(x,\xi,\tau_\epsilon)+\sum_{k\geq 1}\epsilon^kQ_k(\tau,\dots, \partial_x^l\partial_\xi^m\tau_\epsilon)=0,
$ that can be solved recursively under the condition $\partial_\tau h(x,\xi,\tau)|_{h(x,\xi,\tau)=0}\neq 0$.
\section*{Acknowledgements}
We  would like to thank Agis Athanassoulis and Bach-Lien Hua for helpful discussions.
%\centerline{\includegraphics{rosby3.ps}}


\begin{thebibliography}{00}
% please try to use the bibitem system -
% the references should be in alphabetical order of authors' names.
% Articles with a single author first, author will 1 co-author next,
% then author with several co-authors;



\bibitem{CGPS} C. Cheverry, I. Gallagher, T. Paul and L. Saint-Raymond, in preparation.
\bibitem{DU} A. Dutrifoy, A. J. Majda and S. Schochet, A Simple Justification of  
the Singular Limit for Equatorial Shallow-Water Dynamics, 
in Communications on Pure and Applied Math. LXI (2008) 0002-0012.
\bibitem{GA} I. Gallagher and L. Saint-Raymond,  Mathematical study of the  
betaplane model: equatorial waves and convergence results. M\'em. Soc.  
Math. Fr. (N.S.). 107 (2006), v+116 pp.
\bibitem{MA} A. Martinez,  An introduction to semiclassical and microlocal analysis, Springer (2002)

\bibitem{PA} T. Paul, \'Echelles de temps pour l'\'evolution quantique \`a petite constante de Planck, S\'eminaire X-EDP 2007-2008, \'Ecole polytechnique.
\bibitem{P1}J. Pedlosky,  Geophysical fluid dynamics, Springer (1979).

\bibitem{P2} J. Pedlosky, Ocean Circulation Theory, Springer (1996).
\end{thebibliography}
\end{document}